\documentclass[11pt]{amsart}
\usepackage[margin=1.15in]{geometry}
\usepackage{amsmath,amssymb,mathtools,microtype}
\usepackage[hidelinks]{hyperref}
\hypersetup{pdftitle={Divisor Folding and Finite-Group Trace Obstructions},pdfauthor={Abhishek Shankar}}

\newtheorem{theorem}{Theorem}[section]
\newtheorem{proposition}[theorem]{Proposition}
\newtheorem{corollary}[theorem]{Corollary}
\newtheorem{lemma}[theorem]{Lemma}
\theoremstyle{remark}
\newtheorem{remark}[theorem]{Remark}
\newcommand{\Z}{\mathbb Z}
\newcommand{\F}{\mathbb F}
\newcommand{\A}{\mathcal A}
\newcommand{\G}{\mathcal G}
\newcommand{\one}{\mathbf 1}
\DeclareMathOperator{\tr}{tr}

\title{Divisor Folding and Finite-Group Trace Obstructions}
\author{Abhishek Shankar}
\address{OpenProblem.ai}
\email{abhishek.shankar@openproblem.ai}
\thanks{ORCID \href{https://orcid.org/0009-0002-6552-3622}{0009-0002-6552-3622}.}
\date{November 2026}
\subjclass[2020]{42A05, 43A25, 11B30, 05E30}
\keywords{Fourier positivity, divisor cores, resolvent traces, partial difference sets}

\begin{document}
\begin{abstract}
For a symmetric integer frequency set of size $n$ with Fourier floor $K$, we study a positive trace weighted by the reciprocal of the regularized Fourier multiplier. If restriction to frequencies divisible by $p$ leaves $n/p+O(K^2)$ frequencies, then for $n\ge K^4$ and $p\ll n/K^3$ the original trace differs by $O(n)$ from $p$ times the core trace at the inherited regularization. The transfer reaches cores of size $\asymp K^3$ and preserves vanishing triple density with diverging normalized triple count. We construct exactly regular finite abelian spectra at these scales for which the linear core trace bound fails by a power. The examples use unions of classical semiprimitive cyclotomic classes and also give obstructions in both reflection parities in odd characteristic. Finally, we bound the overlap of lifts of a core by its smallest positive Fourier height and the symmetry of the lifts. These bounds exclude balanced lifts of the binary examples in a range of cyclic extensions and quantify the asymmetry required in odd prime split extensions.
\end{abstract}
\maketitle

\section{Introduction}

Let $A=-A\subset\Z\setminus\{0\}$ be finite and let
\[
 F(x)=\sum_{a\in A} e^{2\pi i ax},\qquad P(x)=K+F(x)\ge0,\qquad n=|A|.
\]
Write $D_{ab}=\one_A(a-b)$ for $a,b\in A$. The positivity of $P$ gives $G=KI+D\succeq0$. If $E$ is the orthogonal projection onto the characters indexed by $A$, then the regularized compression $W_H=E(P/(H+P))E$ is positive. The quantity studied here is $\tr(GW_H)$ at the scale $H=(n+K)/K$. This trace retains the reciprocal weight of the Fourier multiplier while using the arithmetic graph on the frequencies.

Restriction to $A\cap p\Z$ averages $P$ over rational translates. It therefore preserves the floor parameter $K$. When this restriction is balanced, the original and restricted traces can be compared in both directions even when the folded triple law has a substantial missing mass. The proof integrates the missing mass against its exact one-coordinate marginal. At the endpoint $p\asymp n/K^3$, it suffices to study a trace of order the core size $K^3$ at regularization of order $K^2$.

The integer setting is motivated by Chowla's cosine problem
\cite{Bourgain,Ruzsa,Sanders,Bedert,Shankar,Quarter}. The quarter-power
bound \cite{Quarter} restricts the integer range \(n\ge K^4\) to
\(K^4\le n\ll K^4\log(2K)\).
The same positivity and balance identities hold in finite abelian groups.
The counterexamples below concern finite groups. Finite abelian spectra with the resulting core scales can have trace much larger than their size. Our examples are made from exactly regular components and have vanishing global triple density. The finite construction underlying one component is the classical semiprimitive cyclotomic construction \cite{CalderbankKantor}, recorded in \cite[Theorem~2.1]{Polhill}. We use its character values to quantify the trace of the union of independent coordinate copies. We then ask which of these examples could arise as cores in a balanced extension. A bound in terms of the smallest positive height constrains the lift whenever enough of its character fibres have repeated heights. For related quotient methods in the study of partial difference sets, see \cite{Zhou}.

Section~\ref{sec:fold} gives the trace transfer and its endpoint consequences. Section~\ref{sec:finite} constructs finite trace obstructions. Section~\ref{sec:lifts} proves the restrictions on cyclic lifts. The results apply whenever the stated balance and symmetry hypotheses hold. The existence of a suitable divisor for every integer spectrum, and the linear trace bound for genuine inherited integer cores, remain open.

\section{Folded traces}\label{sec:fold}

\subsection{Fourier setup}
We work on the compact phase group $X=\widehat{\G}$ with normalized Haar measure, where $\G=\Z$ or $\G$ is a finite abelian group. For $a\in\G$ write $e_a(x)$ for the associated character. Let $A=-A\subset\G\setminus\{0\}$ be finite and nonempty, and let $P=K+F\ge0$, where $F=\sum_{a\in A}e_a$ and $K\ge1$. Set
\[
 n=|A|,\qquad c=K(K-1),\qquad M=n+K,\qquad Z=n+K^3.
\]
Convolution uses normalized Haar measure. Parseval and the exact Fourier coefficients give
\begin{equation}\label{eq:basic}
 \begin{gathered}
  P*P=P+c,\qquad \int P=K,\qquad \int P^2=n+K^2,\\
  T(A):=\int F^3=\#\{(a,b,d)\in A^3:a+b+d=0\}.
 \end{gathered}
\end{equation}
In particular $\rho=P*F^2$ is nonnegative. For $a\in A$ write $r_a=|A\cap(A+a)|$. Then
\[
 \rho=Kn+\sum_{a\in A}r_a e_a.
\]
For $v\ge0$ and $h>0$ put $w_h(v)=v/(h+v)$. Write
\[
 S_A(H)=\int_X\rho(x) w_H(P(x))\,dx.
\]
In the character basis of $E L^2(X)$, compression of multiplication by $P$ has matrix $G=KI+D$. The sum of the entries $G_{ab}$ along the difference $a-b=j$ is $Kn$ for $j=0$, is $r_j$ for $j\in A$, and vanishes otherwise. It follows that $S_A(H)=\tr(GW_H)$. In particular this is a nonnegative trace. The same notation for the empty spectrum denotes the zero trace.

Let $\Gamma\le\G$ have index $p\ge2$, and put $B=A\cap\Gamma$, $m=|B|$. For $\G=\Z$ we take $\Gamma=p\Z$, so $B/p$ is the divided core. The averaging projection $\A_p$ is the average over the annihilator of $\Gamma$ in $X$. Put $P_p=\A_p P=K+F_p$, $F_p=\sum_{b\in B}e_b$, and $M_p=m+K$. We view $P_p$ either on $X$ or as a multiplier on $\widehat\Gamma$, according to the integral. In the integer case, passing to the quotient means writing $P_p(x)$ as a function of $px$.

Let $y_1,\ldots,y_p$ denote the values of $P$ in a single fibre of the restriction $X\to\widehat\Gamma$. Write $Y=pP_p$ and set
\begin{equation}\label{eq:overlap}
 E_p=pP_p^2-\A_p(P^2)=p^{-1}\sum_{i\ne j}y_i y_j,\qquad
 J_p=\int E_p=pm-n+(p-1)K^2.
\end{equation}
These quantities are nonnegative. Also $P\le Y$ and $\int P(Y-P)=J_p$ by averaging the product $P Y$.

\subsection{Triple measures under restriction}
Let $\mu_A$ be the unnormalized measure with density $P(x)P(y)P(x-y)$ on $X^2$. The marginal in each of the variables $x,y,x-y$ has density $P(P+c)$, by \eqref{eq:basic}. In particular $\mu_A(X^2)=Z$. Denote the analogous core measure by $\mu_B$.

\begin{lemma}\label{lem:foldlaw}
The pushforward $\mu_A^{\rm fold}$ of $\mu_A$ to $\widehat\Gamma^2$ satisfies $\mu_A^{\rm fold}\le p\mu_B$. The marginal of $p\mu_B-\mu_A^{\rm fold}$ in the quotient difference variable has density $E_p+(p-1)cP_p$.
\end{lemma}
\begin{proof}
For each of the three quotient variables $u,v,u-v$, divide the lift values by their sum, obtaining probability vectors $\alpha(u)$, $\alpha(v)$, $\alpha(u-v)$. If one sum is zero, both triple densities vanish there. Otherwise averaging the parent product over the $p^2$ pairs of lifts gives $p\Lambda(u,v)$ times the product of the three core heights. Here $\Lambda(u,v)$ is the expectation of the third probability at the difference of the first two lifts, drawn independently according to the first two probability vectors. The choice of fibre representatives may change the difference index by a translation within the annihilator. In every case $0\le\Lambda\le1$. For the marginal, the parent marginal pushes forward to $\A_p(P^2)+cP_p$, whereas $p\mu_B$ has marginal $pP_p^2+pcP_p$.
\end{proof}

\subsection{Two-sided and one-sided comparisons}
\begin{theorem}[Balanced trace transfer]\label{thm:transfer}
For $H>0$ with $H\ge c$,
\begin{equation}\label{eq:transfer}
 \left|S_A(H)-pS_B(H/p)\right|
 \ll pK^3+\frac{M+c}{H}J_p
       +\frac{p}{H}\left(M_p J_p+(p-1)c(m+K^2)\right).
\end{equation}
In particular, fix constants $C_0,C_1$. If $n\ge K^4$, $H=(n+K)/K$, $|m-n/p|\le C_0 K^2$, and $p\le C_1 n/K^3$, then the right side is $O_{C_0,C_1}(n)$. For $p\asymp n/K^3$ and sufficiently large $K$ this compares $O(n)$ parent trace with $O(m)$ core trace, with $m\asymp K^3$ and $H/p\asymp K^2$.
\end{theorem}
\begin{proof}
Replacing $F^2$ by $P^2$ in $S_A(H)$ costs at most
\[
 \int P*|P^2-F^2|\le K\int(2KP+K^2)=3K^3,
\]
since $0\le w_H\le1$ and $P^2-F^2=2KP-K^2$. The cost for $pS_B(H/p)$ is at most $3pK^3$. In the resulting parent trace, write $z=x-y$ and interchange $y$ and $z$. It becomes the integral against $\mu_A$ of $R_o=P(z)/(H+P(x))$. The resulting core trace has test
\[
 R_f=\frac{Y(z)}{H+Y(x)}=\frac{P_p(z)}{H/p+P_p(x)}.
\]
Where the parent triple density is positive,
\begin{equation}\label{eq:ratio}
 R_o-R_f
 =-\bigl(1-P(z)/Y(z)\bigr)R_f
       +R_o\frac{Y(x)-P(x)}{H+Y(x)} .
\end{equation}
The absolute integral of the first term is at most $(M+c)J_p/H$, since $R_f\le Y(z)/H$ and the $z$ marginal is $P(P+c)$. For the second term, integrate in $y$ first. The resulting integral is
\[
 \int_X \frac{P(x)(Y-P)(x)(P*P^2)(x)}{(H+P(x))(H+Y(x))}\,dx
 \le\frac{M}{H}J_p,
\]
using $P*P^2\le M(P*P)=M(P+c)$ and $H\ge c$. Finally, $R_f$ descends to the quotient. Integrating it against the positive difference measure from Lemma~\ref{lem:foldlaw} costs at most
\[
 \frac{p}{H}\int_{\widehat\Gamma}P_p\bigl(E_p+(p-1)cP_p\bigr)
 \le\frac{p}{H}\bigl(M_p J_p+(p-1)c(m+K^2)\bigr).
\]
This proves \eqref{eq:transfer}. In the balanced range $J_p=O(pK^2)$, $M_p=O(n/p)$, and the right side is $O(pK^3+p^2K^5/n)=O(n)$. At the stated endpoint $pm=n+O(n/K)$ and $H/p\asymp K^2$.
\end{proof}

\begin{proposition}[Signed upper comparison]\label{prop:onesided}
For every $H>0$, without the balance hypothesis,
\begin{equation}\label{eq:onesided}
 S_A(H)\le p S_B(H/p)+(p-1)\left(K^3\min(1,pK/H)+2K^2\right).
\end{equation}
\end{proposition}
\begin{proof}
Let $w'=w_{H/p}(P_p)$ and let $\Pi_B$ denote the orthogonal projection onto frequencies in $B$. Positivity and monotonicity give $S_A(H)\le\int (\A_p\rho)w'$. Write $\rho_B=P_p*F_p^2$. At a nonzero frequency $b\in\Gamma$, the coefficient of $E_p$ is
\[
 p\widehat{F_p^2}(b)-\widehat{F^2}(b)+2(p-1)K\one_B(b).
\]
Together with the constant coefficient of $\rho$, this gives
\[
 \A_p\rho-p\rho_B=K(n-pm)+2(p-1)K F_p-\Pi_B E_p.
\]
Since $\Pi_B w'=P_p*w'-K\int w'$ and $E_p,P_p*w'\ge0$, the last term contributes at most $K J_p\int w'$. Now $K(n-pm+J_p)=(p-1)K^3$, $\int w'\le\min(1,pK/H)$, and $\int F_p w'\le K$.
\end{proof}

\subsection{The inherited triple scales}
\begin{corollary}\label{cor:triples}
Under the balanced hypotheses of Theorem~\ref{thm:transfer}, for sufficiently large $K$,
\[
 |p^2T(B)-T(A)|\ll pnK^2,\qquad
 \frac{KT(B)}{m^2}=\rho_*+O\!\left(\rho_*\frac{pK^2}{n}+\frac{pK^3}{n}\right),
 \qquad \rho_*:=\frac{KT(A)}{n^2}.
\]
In particular, at $p\asymp n/K^3$, the conditions $T(A)/n^2\to0$ and $KT(A)/n^2\to\infty$ pass to the core.
\end{corollary}
\begin{proof}
For nonnegative $y_i$ with sum $Y$, expansion gives
\[
 0\le Y^3-\sum_i y_i^3\le\tfrac32Y\bigl(Y^2-\sum_i y_i^2\bigr).
\]
Indeed the coefficient of each $y_i^2y_j$ on both sides of the upper bound is $3$ and the coefficient of each $y_i y_j y_k$ for three distinct indices is $6$ on the left and $9$ on the right. Averaging and using \eqref{eq:overlap},
\[
 0\le p^2\int P_p^3-\int P^3\le\tfrac32 p M_p J_p=O(pnK^2).
\]
Expand $\int P^3=T(A)+3Kn+K^3$, and analogously for $P_p$. The other terms cost $O(pnK^2)$ since $pm=O(n)$ and $pK\ll n$. The second claim follows from $pm=n+O(pK^2)$. At the endpoint its error is $O(\rho_*/K+1)$.
\end{proof}

\begin{remark}
The measure inequality in Lemma~\ref{lem:foldlaw} alone does not give a vanishing bound for the total variation of the normalized laws at the balanced endpoint. The mass fraction of the normalized core law accounted for by this mixture is $a=Z/(p(m+K^3))$. At the principal character $n+K\le p(m+K)$, so
\[
 1-a\ge\frac{(p-1)(K^3-K)}{p(m+K^3)}.
\]
For $m\asymp K^3$ this lower bound stays positive. This concerns the mixture coefficient. It does not give a lower bound for the total variation itself.
\end{remark}

\subsection{A bounded test for the inherited height law}
An envelope $P\le sP_p$ with $s\ge1$ also transfers a bounded logarithmic height test. This comparison can use $s<p$ when a stronger pointwise domination is available.

\begin{proposition}\label{prop:contrast}
Suppose $P\le sP_p$, $H>0$, and $H\ge c$. Put $J_s=\int P(sP_p-P)=s(m+K^2)-(n+K^2)$. For $L>0$ and $a\ge0$ set
\[
 \phi(v)=\min\bigl(\log(1+v/H),\log(1+a/H)\bigr),\qquad
 b_L(t)=\min(1,|t|/L).
\]
Let $\nu=\mu_A/Z$ and let $\nu_{\rm inh}$ be its pushforward to $\widehat\Gamma^2$. Then
\begin{align*}
 \big|&\mathbb E_\nu b_L(\phi(P(x))-\phi(P(y)))
     -\mathbb E_{\nu_{\rm inh}} b_L(\phi(sP_p(u))-\phi(sP_p(v)))\big|\\[-1ex]
 &\hspace{40mm}\le 2L^{-1}\sqrt{J_s/Z} .
\end{align*}
Moreover $\nu_{\rm inh}\le(s/Z)\mu_B$. Its first marginal has density $(\A_p(P^2)+cP_p)/Z$.
\end{proposition}
\begin{proof}
For $t\ge0$, $\log^2(1+t)\le t$, as follows by comparing the derivatives of $\log(1+t)$ and $\sqrt t$. With $Y_s=sP_p$, clipping contracts the difference of the two logarithms, and
\[
 (\phi(Y_s)-\phi(P))^2\le\frac{Y_s-P}{H+P}.
\]
Its integral against $P(P+c)$ is at most $J_s$. Apply Cauchy--Schwarz to each endpoint and the $1/L$ Lipschitz bound of $b_L$. Each lift probability in Lemma~\ref{lem:foldlaw} is at most $s/p$, so $p\Lambda\le s$. The stated marginal follows by pushforward.
\end{proof}

\section{Exactly regular finite examples}\label{sec:finite}

\begin{theorem}[Finite trace obstruction]\label{thm:finite}
Fix $0<\beta<1/3$. There are symmetric sets $A$ in finite abelian groups with true floor $K\to\infty$, such that the graph $D$ is regular and
\[
 m:=|A|\asymp K^3,\qquad \frac{T(A)}{m^2}\longrightarrow0,\qquad
 \frac{KT(A)}{m^2}\longrightarrow\infty.
\]
For $h\asymp K^2$ they satisfy $S_A(h)\gg m q^{2\beta}$, where $q$ is a prime power tending to infinity and $K\asymp q^{(3-\beta)/2}$. The minimum positive value of $P$ is $q^2$. These conclusions hold for powers of two and for odd primes $q$.
\end{theorem}

\subsection{A cyclotomic component}
The single component comes from the semiprimitive construction \cite[Theorem~2.1]{Polhill}. We record the character calculation to specify both the positive height and its multiplicity.

Let $E_0=\F_{q^2}$ and $L_0'=\F_{q^4}$, and let $\chi$ be a multiplicative character on $E_0^\times$ of order $q+1$. Its restriction to $\F_q^\times$ is trivial. Let $\xi=\chi\circ N_{L_0'/E_0}$. Choose a nontrivial additive character of $\F_q$, and extend it to $E_0,L_0'$ using the trace to $\F_q$.

For $1\le j\le q$, the Gauss sum of $\chi^j$ over $E_0$ is $q\chi^j(\omega)$, where $\omega$ is a nonzero element of the trace-zero line over $\F_q$. Indeed the sum of $\chi^j$ on that line minus zero is $(q-1)\chi^j(\omega)$. Its sum on every nonzero trace fibre is the same by scaling, and is $-\chi^j(\omega)$ because the total multiplicative character sum vanishes. Weighting by the additive character of the trace gives the stated Gauss sum. If the characteristic is two we may take $\omega=1$. In odd characteristic $\omega^q=-\omega$, so $\omega^2\in\F_q^\times$. In either case the square of the Gauss sum is $q^2$.

The lifted Gauss sum of $\xi^j$ on $L_0'$ is $-q^2$. To see this, fix trace $s\in E_0$ and norm $t\in E_0^\times$. The number of elements of $L_0'$ with that trace and norm plus the number of ordered pairs in $(E_0^\times)^2$ with sum $s$ and product $t$ equals two. For an irreducible quadratic the two conjugate roots contribute to the first count. For a split quadratic with distinct roots they contribute to the second. For a repeated root the two counts are one each. Weight this identity by the additive character at $s$ and by $\chi^j(t)$, and sum over $s,t$. The contribution from the constant two vanishes and the contribution from the ordered pairs is the square of the Gauss sum over $E_0$.

Let $C$ be the union of $b$ cyclotomic classes of $\xi$, where $1\le b<q+1$, and set
\[
 t=b(q-1),\qquad k=|C|=t(q^2+1),\qquad
 A_0=L_0'\setminus(C\cup\{0\}),\qquad m_0=|A_0|.
\]
Since $\xi(-1)=1$, both $C$ and $A_0$ are symmetric. Identify the characters of the additive group of $L_0'$ with $x\in L_0'$, so that the character at $x$ has value at $a$ given by the chosen additive character of $ax$. Expanding the class indicators in powers of $\xi$, the principal multiplicative term contributes $-b/(q+1)$ at a nonzero $x$. The other terms give $-q^2(\one_C(x^{-1})-b/(q+1))$. Thus
\begin{equation}\label{eq:twoheights}
 F_C(x)=t-q^2\one_C(x^{-1})\qquad (x\ne0).
\end{equation}
For $K_0=t+1$, $P_0=K_0+F_{A_0}$ has value $q^2$ at exactly $k$ nonprincipal phases and vanishes at the other nonprincipal phases. Its principal value is $K_0+m_0$. Consequently $K_0$ is the true floor and the smallest positive height is $q^2$.

For each $a\in A_0$, the difference representation count in $C$ is $t(t+1)$. In fact $(F_C-t)(F_C-t+q^2)$ vanishes off the principal character by \eqref{eq:twoheights}, and at the principal character it equals
\[
 (k-t)(k-t+q^2)=q^4t(t+1).
\]
Taking its Fourier coefficient at $a\in A_0$ gives the claim since $a\notin C\cup\{0\}$. Complementing both entries in the difference count, the degree of each vertex in $A_0$ is
\begin{equation}\label{eq:degree}
 d_0=q^4-2(k+1)+t(t+1).
\end{equation}

\subsection{Independent copies and the trace}
\begin{proof}[Proof of Theorem~\ref{thm:finite}]
Take $b=\lfloor q^\beta\rfloor$ and $\ell=\lfloor q^{(1-3\beta)/2}\rfloor$. Use $\ell$ copies of $A_0$ in independent coordinate axes of $(L_0')^\ell$. There are no mixed triples. In the phase group, $P=\sum_{i=1}^\ell P_0(x_i)$ has true floor $K=\ell K_0$ since the summands can vanish simultaneously. The smallest positive height is $q^2$ since one coordinate can have that height while the others vanish. The graph is the disjoint union of the component graphs, so it is regular with degree $d_0$ and $T(A)=md_0$. Equations \eqref{eq:twoheights}--\eqref{eq:degree} give
\[
 m_0\asymp d_0\asymp q^4,\quad m=\ell m_0\asymp K^3,
 \quad T(A)/m^2\asymp\ell^{-1},\quad KT(A)/m^2\asymp K_0.
\]
For a regular graph $\rho=K(m-d_0)+d_0P$. The first term is nonnegative. Since $w_h$ is increasing and $P\ge P_0(x_i)$,
\[
 S_A(h)\ge d_0\sum_{i=1}^{\ell}\int P_0(x_i)w_h(P_0(x_i))
           =d_0\ell\int P_0w_h(P_0).
\]
The medium heights alone give $\int P_0w_h(P_0)\ge k/(h+q^2)\gg q^{2\beta}$ since $k\asymp q^{3+\beta}$ and $h\asymp q^{3-\beta}$. This proves the lower bound.
\end{proof}

\begin{remark}\label{rem:integer}
The dense components of this construction use torsion. For a symmetric set of $2s$ nonzero integers, every ordered positive relation $a+b=d$ is determined by the positive sum $d$ and the first summand $a<d$. There are at most $s(s-1)/2$ such relations, and every ordered zero-sum triple arises from one of them with a factor of six. Thus its triple count is at most $3s(s-1)$. In the components above $T(A_0)/m_0^2=d_0/m_0\to1$. The global sparse triple conditions alone do not capture this local restriction.
\end{remark}

\subsection{Both reflection parities in odd characteristic}
On the character space indexed by $A$, let $\mathcal R$ send a vector $(u_a)$ to $(u_{-a})$ and let $\Pi_\pm=(I\pm\mathcal R)/2$. Both $G$ and $W_h$ commute with reflection since $A=-A$ and $P$ is even.

\begin{proposition}\label{prop:parity}
For the examples of Theorem~\ref{thm:finite} along odd primes $q$, each of the two traces on the reflection eigenspaces satisfies
\[
 \tr\bigl(\Pi_\pm G W_h\Pi_\pm\bigr)\gg m q^{2\beta}.
\]
\end{proposition}
\begin{proof}
In one component put $C_0=C\cup\{0\}$ and write $v_x=(e_a(x))_{a\in A_0}$. Fix a medium-height phase $x$, so $x\ne0$ and $P_0(x)=q^2$. By expanding the two complement indicators in the adjacency product, for $a\in A_0$,
\[
 (D_{A_0}v_x)_a=-F_{C_0}(x)(1+e_a(x))
    +\sum_{b\in L_0'}\one_{C_0}(a-b)\one_{C_0}(b)e_b(x).
\]
Here $-F_{C_0}(x)=q^2-O(q^{1+\beta})$. By Parseval on the frequency group with counting measure, the squared norm of the last term on the whole group is
\[
 q^{-4}\sum_{z\in L_0'} |F_{C_0}(z)F_{C_0}(z+x)|^2.
\]
The principal exceptions $z=0,-x$ contribute $O(q^{6+2\beta})$. The terms with at least one medium nonprincipal value contribute $O(q^{7+\beta})$, since there are at most $2k$ of them. The remaining terms contribute $O(q^{4+4\beta})$. The normalized sum is $o(q^8)$ uniformly on these phases, and thus
\[
 D_{A_0}v_x=q^2(\one+v_x)+o(q^4)
\]
in Euclidean norm. Reflection commutes with $D_{A_0}$, and
\[
 \|\Pi_\pm v_x\|^2=\tfrac12\bigl(m_0\pm F_{A_0}(2x)\bigr)
                       =\tfrac12m_0+O(q^2)
\]
since $2x\ne0$. In the even space the inner product of $\one$ with $\Pi_+v_x$ is $F_{A_0}(x)=O(q^2)$. In the odd space the constant term vanishes. It follows in either parity that the quadratic form of $D_{A_0}$ at $\Pi_\pm v_x$ is $\gg q^6$, since the norm of $v_x$ is $O(q^2)$ and the error contributes $o(q^6)$. Adding $K I$ preserves the lower bound.

The global evaluation vector has one component of this form for each coordinate. Its parity quadratic form under $G$ is the sum of the corresponding component forms, all nonnegative by positivity of $G$. Integrate each component over its own medium-height phases and use $w_h(P)\ge w_h(q^2)$ there. The contribution per component is at least a constant times $(k/q^4)q^6 w_h(q^2)\gg q^4 q^{2\beta}$. Summing proves the result.
\end{proof}

\section{Lift restrictions from a positive height gap}\label{sec:lifts}

Let $\G$ now be finite, let $\Gamma\le\G$ have index $p\ge2$, and let $B=A\cap\Gamma$ be a nonempty core. Write $\kappa=-\min F_B$ for its true floor and let $L_0$ be the smallest positive value of $\kappa+F_B$ on $\widehat\Gamma$. Let the parent use any $K'\ge\kappa$ such that $P=K'+F_A\ge0$. Its lift-overlap $J_p(P)$ is defined by \eqref{eq:overlap} with $K'$ in place of $K$.

\begin{theorem}[Repeated heights in a fibre]\label{thm:symmetry}
Let $\mathcal H$ be a group of automorphisms of $\G$ fixing $\Gamma$ pointwise. Let $\bar P$ be the orthogonal average of $P$ under the induced action on $\widehat\G$, and put
\[
 U=\sum_{\sigma\in\mathcal H}(\sigma-\mathrm{id})\G.
\]
If every nonfixed character orbit has size at least $r\ge2$, then
\begin{equation}\label{eq:symmetry}
 J_p(P)+\|P-\bar P\|_2^2
 \ge(1-1/r)L_0\left[p\kappa-\frac{p}{|U|}\bigl(K'+|A\cap U|\bigr)\right]_+.
\end{equation}
\end{theorem}
\begin{proof}
The average $\bar P$ is nonnegative, has the same mean and the same core multiplier as $P$, and is constant on each orbit within a lift fibre. Consequently
\[
 J_p(\bar P)=p\int P_p^2-\int\bar P^2=J_p(P)+\|P-\bar P\|_2^2.
\]
Let $y_j$ now denote the lift values of $\bar P$. On a nonfixed orbit one has $y_j\le Y/r$ since the other orbit values agree. In the formula $E_p=p^{-1}\sum_j y_j(Y-y_j)$, the contribution from these lifts is at least $(1-1/r)P_p$ times their total height.

The fixed characters are exactly $U^\perp$. The total integrated height on the fixed lifts is
\[
 \frac{p}{|\G|}\sum_{\psi\in U^\perp}\bar P(\psi)
 =\frac{p}{|U|}\bigl(K'+|A\cap U|\bigr).
\]
For the equality, $\bar P=P$ at every fixed character, and orthogonality makes the sum of $e_a$ over $U^\perp$ vanish unless $a\in U$. The integrated height of all lifts is $pK'$. On phases where $\kappa+F_B$ vanishes, their total height is $p(K'-\kappa)$ per fibre. Removing these phases and the fixed lifts therefore leaves integrated nonfixed height at least the positive part in \eqref{eq:symmetry}. On the remaining phases $P_p\ge L_0$, so integration of the lower bound for $E_p$ proves the result.
\end{proof}

\begin{corollary}[Binary cores in cyclic extensions]\label{cor:binary}
If $\Gamma$ has exponent two and $\G/\Gamma$ is cyclic of order $p$, then
\begin{equation}\label{eq:binary}
 J_p(P)\ge\frac{L_0}{2}\bigl[p\kappa-2(K'+p-1)\bigr]_+.
\end{equation}
For the binary examples of Theorem~\ref{thm:finite}, fix $C\ge1$. If $K'\le C\kappa$ and $p\ge8C$, then as $q\to\infty$ this is incompatible with $J_p(P)=O(p(K')^2)$ with a uniform constant.
\end{corollary}
\begin{proof}
Use inversion as the automorphism. It fixes $\Gamma$ pointwise, and symmetry of $A$ gives $P=\bar P$. The nonfixed orbits have size two and $U=2\G$. If $g$ lifts a quotient generator, then $2\G=\langle2g\rangle$ since $2\Gamma=0$. Its order lies between $p/2$ and $p$ because its quotient image has order $p/\gcd(2,p)$ and $2pg=0$. Hence $(p/|U|)(K'+|A\cap U|)\le2(K'+p-1)$. For the examples, $L_0=q^2$ and $\kappa\asymp q^{(3-\beta)/2}=o(L_0)$. The bracket in \eqref{eq:binary} is at least a fixed multiple of $p\kappa$ for sufficiently large $q$. This yields $J_p(P)\gg p\kappa L_0$.
\end{proof}

\begin{corollary}[Odd prime split extensions]\label{cor:split}
Suppose $\G=\Gamma\times\Z/p\Z$ for an odd prime $p$. Let $\bar P$ average over all nonzero multipliers of the quotient coordinate. Then
\begin{equation}\label{eq:split}
 J_p(P)+\|P-\bar P\|_2^2\ge
 \frac{p-2}{p-1}L_0[p\kappa-K'-(p-1)]_+ .
\end{equation}
If $q_g=\#\{t\ne0:(g,t)\in A\}$ for $g\in\Gamma$, then
\begin{equation}\label{eq:variance}
 \|P-\bar P\|_2^2=\sum_{g\in\Gamma}q_g\bigl(1-q_g/(p-1)\bigr).
\end{equation}
In particular, for the examples of Theorem~\ref{thm:finite}, if $K'\le C\kappa$, $p\ge 8C$, and $J_p(P)=O(p(K')^2)$, then the variance in \eqref{eq:variance} is $\gg p\kappa L_0$ for sufficiently large $q$.
\end{corollary}
\begin{proof}
The quotient multipliers have nonfixed character orbits of size $p-1$, and $U=\{0\}\times\Z/p\Z$. Apply \eqref{eq:symmetry} and $|A\cap U|\le p-1$. Averaging in the frequency group leaves the coefficient at $(g,0)$ unchanged and replaces the $p-1$ coefficients at $(g,t)$, $t\ne0$, by their mean $q_g/(p-1)$. Parseval yields \eqref{eq:variance}. In the examples $L_0/\kappa\to\infty$, so the overlap bound is lower order than the right side of \eqref{eq:split} under the stated assumptions.
\end{proof}

The balance condition of Theorem~\ref{thm:transfer} implies $J_p(P)=O(p(K')^2)$ in these finite extensions. Corollary~\ref{cor:binary} therefore excludes such a lift of the binary examples at comparable floors. Corollary~\ref{cor:split} gives a quantitative asymmetry requirement for split lifts and allows noninvariant parents. Inversion of integer phases generally moves the underlying core character, so the fibre argument for binary cores requires its finite-group hypothesis.

\end{document}